\documentclass[11pt,A4paper]{article}
\usepackage{amsmath,amsthm,amssymb,amsfonts}
\usepackage[numbers,sort&compress]{natbib}
\usepackage{enumerate}
\usepackage{color,colordvi}
\usepackage{soul}
\usepackage{fullpage}

\usepackage{changes}

\usepackage[letterpaper,top=2cm,bottom=2cm,left=3cm,right=3cm,marginparwidth=2cm]{geometry}

\usepackage{amsmath,amsfonts,amsthm,mathtools,bm}
\usepackage{graphicx}
\usepackage[colorlinks=true, allcolors=blue]{hyperref}
\usepackage[capitalise]{cleveref}
\usepackage{standalone}
\usepackage{algorithm} 
\usepackage{indentfirst}
\usepackage{enumerate}
\usepackage[affil-it]{authblk}
\usepackage{diagbox}
\usepackage{changes}
\usepackage{enumitem}
\usepackage{algpseudocode}
\usepackage{bbm}
\usepackage{orcidlink}

\algrenewcommand\algorithmicrequire{\textbf{Input:}}
\algrenewcommand\algorithmicensure{\textbf{Output:}}

\usepackage[mathscr]{euscript}
\usetikzlibrary{calc}

\usetikzlibrary{decorations.pathreplacing} 

\theoremstyle{plain}
\newtheorem{lemma}{Lemma}[section]
\newtheorem{theorem}[lemma]{Theorem}

\theoremstyle{definition}

\newtheorem{case}{Case}
\newtheorem{subcase}{Subcase}[case]

\newtheorem{claim}{Claim}

\newtheorem{definition}[lemma]{Definition}

\newtheorem{problem}[lemma]{Problem}

\newtheorem{remark}[lemma]{Remark}

\newtheorem*{claim*}{Claim}

\numberwithin{equation}{section}

\newcommand{\RR}{\mathbb{R}}

\newcommand{\RomanNum}[1]{\uppercase\expandafter{\romannumeral #1}}

\DeclarePairedDelimiter{\card}{\lvert}{\rvert}

\title{A Faber--Krahn inequality for trees}

\author{Huiqiu Lin\footnote{email: huiqiulin@126.com}~\orcidlink{0000-0002-6072-4647}}
\author{Lianping Liu\footnote{email: y12252216@mail.ecust.edu.cn}~\orcidlink{0009-0007-9378-5028}}
\author{Zhe You\footnote{email: y30231280@mail.ecust.edu.cn}~\orcidlink{0009-0006-9769-5372}}
\affil{School of Mathematics, East China University of Science and Technology, 130 Meilong Road, Shanghai 200237, China.}

\date{}

\begin{document}
\maketitle

\begin{abstract}
The well-known Faber-Krahn theorem states that the ball has the lowest first Dirichlet eigenvalue among all domains of the same volume in $\mathbb{R}^n$. 
Leydold (Geom. Funct. Anal, 1997) gave the discrete version of Faber-Krahn inequality for regular trees with boundary.
B{\i}y{\i}ko{\u{g}}lu and Leydold (J. Combin.
Theory Ser. B, 2007) demonstrated that the Faber--Krahn inequality holds for the class of trees with boundary with the same degree sequence. They further posed the following question: Give a characterization of all graphs in a given class \(\mathcal{C}\) with the Faber-Krahn property.
In this paper, we show the Faber-Krahn property for trees with given matching number. 
Our result can imply the Klob\"ur\v{s}tel theorem, i.e., the Faber-Krahn inequality for trees with given number of interior vertices and boundary vertices.
\end{abstract}

\section{Introduction}
Dirichlet eigenvalue problem constitutes a fundamental topic within spectral geometry. 
For a bounded domain $\Omega$ with smooth boundary $\partial \Omega$ in a Riemannian manifold, the Dirichlet eigenvalue problem is shown as
\begin{align*}
    \Delta f + \lambda f = 0 \text{ in $\Omega$}
\end{align*}
with Dirichlet boundary condition 
\begin{align*}
    f|_{\partial\Omega} = 0.
\end{align*}
This eigenvalue problem can also be viewed as the eigenvalue problem of the Dirichlet Laplacian operator $\Delta_\Omega$, which is defined by $\Delta_\Omega f \coloneqq (\Delta \widehat{f})|_\Omega$, where $\widehat{f}$ is an extension of $f$ to the whole manifold by assigning 0 to the region outside $\Omega$. 

Let $G = (V, E)$ be an undirected simple graph, where $V, E$ are respectively the set of vertices and edges of $G$. 
A graph with boundary is a pair $(G, B)$, where $B \subsetneq V$ is a non-empty set. 
The set $B$ is called the boundary and $\Omega \coloneqq V \setminus B$ is called  the interior. 
If $G$ is a tree, then we denote the pair $(G,B)$ by $T$ for simplicity, where we choose all leaves (vertices of degree one) as boundary vertices.
The Dirichlet eigenvalue problem on a graph $G$ with boundary $B$ is defined as follows.
\begin{equation*}
    \begin{cases}
        -\Delta f(x) = \lambda f(x), & x\in \Omega, \\ 
         f(x) = 0, & x\in B,
    \end{cases}
\end{equation*}
where $\Delta f(x) = \sum\limits_{y \sim x} (f(y)-f(x))$.
We call $\lambda$ the Dirichlet eigenvalue and $f$ the eigenfunction corresponding to $\lambda$.
The Dirichlet eigenvalues are ordered by $0 < \lambda_{1} \leq \lambda_{2} \leq \ldots \leq \lambda_{|\Omega|}$. 
For any $0 \neq f \in \RR^{\Omega}$, the Rayleigh quotient of Dirichlet Laplacian is defined as
\begin{align*}
    R_{(G,B)}(f) &\coloneqq \frac{\sum\limits_{(x, y) \in E(G) }(\widehat{f}(x)-\widehat{f}(y))^{2}}{\sum\limits_{x \in \Omega} f^{2}(x)},
\end{align*}
where $\widehat{f}$ is obtained by extending $f$ to $\RR^V$ by $0$.

The variational characterizations of $\lambda_{k}$ are given by
\begin{align*}
    \lambda_{k}(G,B) &= \min_{\substack{W \subset \mathbb{R}^{\Omega}\\ \dim W=k}} \max_{\substack{0 \neq f \in W}} R_{(G,B)}(f).
\end{align*}

In the past few decades, significant progress has been made in extending classical eigenvalue problems from continuous to discrete settings. 
Various analogues of inequalities for Dirichlet and Neumann Laplacian eigenvalues have been established. 
For instance, on planar domains, Rohleder~\cite{MR4958511} obtained fundamental inequalities relating Neumann and Dirichlet eigenvalues, thereby enriching the classical spectral theory for bounded domains.
Beyond these specific extensions, comparisons among Dirichlet, Neumann, and Laplacian eigenvalues on graphs have drawn increasing interest. 
Shi and Yu~\cite{MR4936329} conducted a systematic study of such eigenvalue comparisons and further examined their rigidity, offering fundamental insights into spectral analysis on graphs~\cite{MR4974518}.
Many scholars have investigated corresponding problems on graphs. 
For induced subgraphs of the $n$-dimensional integer lattice, Bauer and Lippner~\cite{bauer2022eigenvalue} considered a discrete version of Pólya's conjecture. Furthermore, Hua and Li~\cite{MR4956579} extended Bauer and Lippner's result to poly-Laplace operators.
Hua, Lin and Su~\cite{hua2023payne} proved some analogues of Payne–Pólya-Weinberger, Hile–Protter and Yang's inequalities. For Cayley graphs of finitely generated amenable groups, Hua and Yadin~\cite{hua2024universal} extended the result of Hua, Lin and Su, and this result also applies to regular trees.
Wang and Hou~\cite{wang2025faber} further advanced this line of research by establishing a Faber-Krahn type inequality for supertrees. 

The study of Faber-Krahn type inequalities constitutes an important problem for various special classes of graphs.
B{\i}y{\i}ko{\u{g}}lu and Leydold~\cite{biyikouglu2007faber} proposed the following problem.
\begin{problem}[{\cite[Problem 1]{biyikouglu2007faber}}]\label{pro:Faber-Krahn}
    Give a characterization of all graphs in a given class $\mathcal{C}$ with the $\mathit{Faber}$–$\mathit{Krahn}$ property, i.e., characterize those graphs in $\mathcal{C}$ which have minimal first Dirichlet eigenvalue for a given “volume”.
\end{problem}
Let $\mathcal{T}^{(n,k)}$ denote the set of trees of order $n$ with $k$ interior vertices.
Klob\"ur\v{s}tel theorem tells us the Faber-Krahn property in the class $\mathcal{T}^{(n,k)}$.

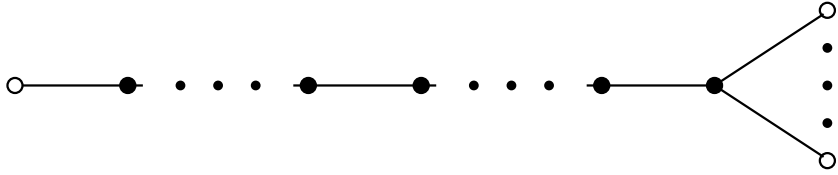
\begin{figure}[htbp]
    \centering
    \begin{tikzpicture}[    
x=1.00mm, y=1.00mm,    
inner xsep=0pt, inner ysep=0pt,    
outer xsep=0pt, outer ysep=0pt,    
node circle/.style={line width=0.30mm, draw=L, fill=F, circle},    
node circle small/.style={line width=0.15mm, draw=L, fill=F, circle}
]
\path[line width=0mm] (41.90,36.47) rectangle +(66.10,26.59);
\definecolor{L}{rgb}{0,0,0}
\definecolor{F}{rgb}{0,0,0}
\path[line width=0.30mm, draw=L] (40,50) circle (1.00mm);
\path[line width=0.30mm, draw=L, fill=F] (55,50) circle (1.00mm);
\path[line width=0.30mm, draw=L, fill=F] (62,50) circle (0.5mm);
\path[line width=0.30mm, draw=L, fill=F] (67,50) circle (0.5mm);
\path[line width=0.30mm, draw=L, fill=F] (72,50) circle (0.5mm);
\path[line width=0.30mm, draw=L, fill=F] (79,50) circle (1.00mm);
\path[line width=0.30mm, draw=L, fill=F] (94,50) circle (1.00mm);
\path[line width=0.30mm, draw=L, fill=F] (101,50) circle (0.5mm);
\path[line width=0.30mm, draw=L, fill=F] (106,50) circle (0.5mm);
\path[line width=0.30mm, draw=L, fill=F] (111,50) circle (0.5mm);
\path[line width=0.30mm, draw=L, fill=F] (118,50) circle (1.00mm);
\path[line width=0.30mm, draw=L, fill=F] (133,50) circle (1.00mm);
\path[line width=0.30mm, draw=L] (148,60) circle (1.00mm);
\path[line width=0.30mm, draw=L, fill=F] (148,55) circle (0.5mm);
\path[line width=0.30mm, draw=L, fill=F] (148,50) circle (0.5mm);
\path[line width=0.30mm, draw=L, fill=F] (148,45) circle (0.5mm);
\path[line width=0.30mm, draw=L] (148,40) circle (1.00mm);
\path[line width=0.30mm, draw=L] (41,50) -- (57,50);
\path[line width=0.30mm, draw=L] (77,50) -- (96,50);
\path[line width=0.30mm, draw=L] (116,50) -- (133,50);
\path[line width=0.30mm, draw=L] (133,50) -- (147.5,59.5);
\path[line width=0.30mm, draw=L] (133,50) -- (147.5,40.5);
\end{tikzpicture}
    \caption{A comet has the Faber-Krahn property in class $\mathcal{T}^{(n,k)}$.}
    \label{a comet}
\end{figure}

\begin{theorem}[Klob\"ur\v{s}tel theorem {\cite[Theorem 1]{biyikouglu2007faber}}]\label{Klobürštel theorem}
    A tree  $T$  has the Faber-Krahn property in the class $\mathcal{T}^{(n,k)}$ if and only if $T$ is a star with a long tail, i.e., a comet, see~\cref{a comet}. 
    $T$ is then uniquely determined up to isomorphism.
\end{theorem}

A matching in a graph is a set of disjoint edges.
The (maximum) matching number of a graph $G$ is the maximum number of edges among all matchings in $G$. 
We investigate Problem~\ref{pro:Faber-Krahn} in the following classes of trees with boundary:
\begin{align*}
    &\mathcal{T}(n,m) = \{T: T \text{ is a tree with } n  \text{ vertices and matching number } m \},\\
    &\mathcal{T}(n,m,b) = \{T \in \mathcal{T}(n,m) : \text{the number of leaves of } T \ \text{is} \ b \}.
\end{align*}

We need the following notions.
Let $y \sim x$ denote that $y$ is adjacent to $x$. 
We write $N_G(x):=\{\,y\in V(G) : xy\in E(G)\,\}$ for the set of neighbors of
$x$ in $G$ and $|\cdot|$ for the cardinality of a set.
Denote the degree of $x$ in $G$ by $d_G(x):=\card{N_G(x)}$.
A path of length $l$ in a graph is a sequence of distinct vertices $v_0,v_1,\cdots, v_l$ such that $v_{i-1} \sim v_i$ for all $i = 1,2,\ldots,l$. 
The (geodesic) distance between two vertices $x,y$, which is the length of a shortest path connecting $x$ and $y$, is denoted by $d(x,y)$. 
The diameter $D$ of a graph $G$ is defined by $D = \sup_{v, w \in V} d(v,w)$. 
For a vertex set $Y$, we write $d(x,Y):=\inf\{d(x,y)|y\in Y\}$. 
The inscribed radius of $(G, B)$ is defined by $r = \sup_{v \in V} d(v, B)$. 
We proceed to introduce the notion of  balls and ball approximations in a graph.
A $ball$ $approximation$  $B(v_0, \{r,r+1\})$  with center $v_0$ and radius $r \in \mathbb{N}$ is a connected graph where every boundary vertex $w$ has geodesic distance $d(v_0, w) \in  \{r,r+1\}$. 
In particular, if every boundary vertex $w$ satisfies $d(v_0, w) = r$, then the set $B(v_0, r)$ is a ball.

Now we consider a family of trees as follows.
\begin{definition}\label{def:Leafed Path Tree}
    Let $p$, $q$, $b$ be integers such that $p \geq 0$ and $b \geq q \geq 2$.
    The tree $T(p, q, b)$ is obtained through the following process. 
    Let $u_{1} \sim \cdots u_{2} \cdots \sim u_{p+1} \sim u_{p+2}\cdots\sim u_{p+q}$ be a path of length $p+q-1$. 
    Attach a vertex to $u_{1},u_{p+2},u_{p+3},\cdots ,u_{p+q-1}$ respectively, and attach $b+1-q$ vertices to $u_{p+q}$.
    (See~\cref{def:Leafed Path Tree})
\end{definition}
It is easy to check that the matching number of $T(p, q, b)$ is $q+ \lfloor\frac{p}{2}\rfloor$.
    
\begin{figure}[htbp]
    \centering
       \begin{tikzpicture}[    
x=1.00mm, y=1.00mm,    
inner xsep=0pt, inner ysep=0pt,    
outer xsep=0pt, outer ysep=0pt,    
node circle/.style={line width=0.30mm, draw=L, fill=F, circle},    
node circle small/.style={line width=0.15mm, draw=L, fill=F, circle}
]
\path[line width=0mm] (41.90,36.47) rectangle +(66.10,26.59);
\definecolor{L}{rgb}{0,0,0}
\definecolor{F}{rgb}{0,0,0}
\path[line width=0.30mm, draw=L] (40,60) circle (1.00mm);
\path[line width=0.30mm, draw=L, fill=F] (40,50) circle (1.00mm);
\path[line width=0.30mm, draw=L, fill=F] (55,50) circle (1.00mm);
\path[line width=0.30mm, draw=L, fill=F] (65,50) circle (0.5mm);
\path[line width=0.30mm, draw=L, fill=F] (70,50) circle (0.5mm);
\path[line width=0.30mm, draw=L, fill=F] (75,50) circle (0.5mm);
\path[line width=0.30mm, draw=L, fill=F] (85,50) circle (1.00mm);
\path[line width=0.30mm, draw=L, fill=F] (100,50) circle (1.00mm);
\path[line width=0.30mm, draw=L, fill=F] (110,50) circle (0.5mm);
\path[line width=0.30mm, draw=L, fill=F] (115,50) circle (0.5mm);
\path[line width=0.30mm, draw=L, fill=F] (120,50) circle (0.5mm);
\path[line width=0.30mm, draw=L, fill=F] (130,50) circle (1.00mm);
\path[line width=0.30mm, draw=L, fill=F] (145,50) circle (1.00mm);
\path[line width=0.30mm, draw=L] (155,60) circle (1.00mm);
\path[line width=0.30mm, draw=L, fill=F] (155,55) circle (0.5mm);
\path[line width=0.30mm, draw=L, fill=F] (155,50) circle (0.5mm);
\path[line width=0.30mm, draw=L, fill=F] (155,45) circle (0.5mm);
\path[line width=0.30mm, draw=L] (155,40) circle (1.00mm);
\path[line width=0.30mm, draw=L] (40,50) -- (40,59);
\path[line width=0.30mm, draw=L] (40,50) -- (57,50);
\path[line width=0.30mm, draw=L] (83,50) -- (102,50);
\path[line width=0.30mm, draw=L] (128,50) -- (145,50);
\path[line width=0.30mm, draw=L] (154.5,59.5) -- (145,50);
\path[line width=0.30mm, draw=L] (154.5,40.5) -- (145,50);
\path[line width=0.30mm, draw=L] (100,50) -- (100,59);
\path[line width=0.30mm, draw=L] (130,50) -- (130,59);
\path[line width=0.30mm, draw=L] (100,60) circle (1.00mm);
\path[line width=0.30mm, draw=L] (130,60) circle (1.00mm);
\node at (40,47){$u_1$};
\node at (55,47){$u_2$};
\node at (85,47){$u_{p+1}$};
\node at (100,47){$u_{p+2}$};
\node at (130,47){$u_{p+q-1}$};
\node at (144,47){$u_{p+q}$};
\node at (161,60){$u_{p+q,1}$};
\node at (165,40){$u_{p+q,b+1-q}$};
\end{tikzpicture}
    \caption{The tree $T(p, q, b)$}
    \label{Leafed Path Tree}
\end{figure}
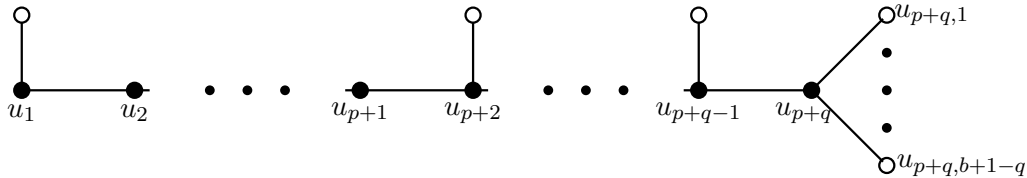

The famous Faber–Krahn theorem gives an important isoperimetric inequality for Dirichlet eigenvalues. 
The theorem states that the ball has the lowest first Dirichlet eigenvalue among all bounded domains with the same volume in $\mathbb{R}^n$ with the standard Euclidean metric. 
It was first conjectured by Rayleigh and proved independently by Faber~\cite{faber1923beweis} and Krahn~\cite{krahn1925rayleigh} for $\mathbb{R}^2$. 
A proof of its generalized version can be found in~\cite{chavel1984eigenvalues}.
In the discrete setting, Friedman~\cite{friedman1993some} conjectured that a Faber–Krahn type inequality should hold for regular trees. 
However, this conjecture was shown to be false, see~\cite{leydold1996faber, pruss1998discrete}. 
Subsequently, Leydold~\cite{leydold2002geometry} established a Faber–Krahn inequality for regular trees and gave a characterization of the extremal trees, which are ball approximations.
B{\i}y{\i}ko{\u{g}}lu and Leydold~\cite{biyikouglu2007faber} also show that trees that have lowest first Dirichlet eigenvalue for a given degree sequence are ball approximations.
Therefore balls are of particular interest for our investigations.
Notice that $T(p,2,b)$ is ball approximation.
Next, we can formulate Faber-Krahn type theorem for the classes $\mathcal{T}(n,m)$.
\begin{theorem}\label{cor:m_extremal tree}
   Let $n\geq 3$ and $m\geq 1$ be integers.
   A tree $T$ has the Faber-Krahn property in the class $\mathcal{T}(n,m)$ 
    if and only if one of the following holds.
        \begin{enumerate}
        \item For $m=1$, $T$ is isomorphic to $T(0,\,1,\,n-1)$.
        \item $m \geq 2$
        \begin{enumerate}
        \item if $n \geq 2m+1$, then $T$ is isomorphic to $T(2m-3,\,2,\,n+1-2m)$.
        \item if $n = 2m$, then $T$ is isomorphic to $T(2m-4,\,2,\, 2)$.
        \end{enumerate}
        \end{enumerate}
\end{theorem}
However, for the classes $\mathcal{T}(n,m,b)$, the trees that satisfy the Faber–Krahn property may be balls, ball approximations, or other structures that are not ball-like at all.
For $m=1$, we have $b=n-1$ and $T$ is isomorphic to $T(0,\,1,\,n-1)$. So in the following, we only need to consider that $m \geq 2$.
\begin{theorem}\label{thm:m_b_extremal tree}
    Let $m\geq2$, $n\geq 3$, $b\geq 2$ be integers and $t = 2m + b - n$.
    If $T \in \mathcal{T}(n,m,b)$, then $1 \leq t\leq \min\{b,m\}$.
    A tree $T$ has the Faber-Krahn property in the class $\mathcal{T}(n,m,b)$ if and only if one of the following holds.
    \begin{enumerate}
        \item For $ t= 1$, $T$ is isomorphic to $T(2m -3,2,b)$.
        \item For $ t \geq 2$,
        \begin{enumerate}
            \item if $t < m$, then $T$ is isomorphic to $T(2m-2t, t,b)$.
            \item if $t = m =b$, then $T$ can be any tree satisfying that each interior vertex is adjacent to exactly one leaf.
            \item if $t = m < b$, then $T$ is isomorphic to $T(0, t,b)$.
            \end{enumerate}
        \end{enumerate}
\end{theorem}

\section{Preliminaries}
In this section, we provide some techniques which are important in our proof.
For the first Dirichlet eigenvalue and its eigenfunction of a general graph, we have following properties. 
\begin{lemma}[{\cite[Lemma 2.1]{lin2025estimatesdirichleteigenvaluegraphs}}]\label{lem:simple eigenvalue-positive}
    Let $G = (V, E)$ be a graph with boundary $B$ and $\Omega = V \setminus B$ be the interior. 
    Suppose that the interior is connected. 
    Let $\lambda_1(G,B)$ be the first Dirichlet eigenvalue. 
    Then the following statements hold.
    \begin{enumerate}
        \item[(i)] The eigenvalue $\lambda_1$ is a simple eigenvalue;
        \item[(ii)] there exists a positive eigenfunction corresponding to $\lambda_1$;
        \item[(iii)] all non-negative eigenfunctions of the Dirichlet Laplacian operator must belong to $\lambda_1$.
    \end{enumerate}
\end{lemma}

We say that $G$ has separated boundary if $B$ has a distinct vertex for each boundary edge, which means the endpoints of it consist of one boundary vertex and one interior vertex.
If $G_1$ and $G_2$ are graphs with separated boundary, we say that $G_2$ is an extension of $G_1$, written $G_1 \subseteq G_2$, if there exists an isometric embedding of the realization of $G_1$ into $G_2$ which preserves the degree of each interior vertex; that is, $G_1$ is obtained from $G_2$ by declaring some of its interior vertices to be boundary vertices and by shortening some of the boundary edges' lengths (and ignoring any edges now connecting two boundary vertices).  
For $G_1 \subseteq G_2$, we recall the following lemma about Dirichlet eigenvalues of $G_1$ and $G_2$.

\begin{theorem}\label{thm:subgraph_larger_original graph}
Given two graphs $G_1$ and $G_2$, if $G_1 \subseteq G_2$, then the $k$-th eigenvalue of $G_2$ for given boundary conditions is
$\le$ that of $G_1$ with the Dirichlet induced boundary conditions; in particular, this holds
for the Dirichlet eigenvalues of $G_1$ and $G_2$.
\end{theorem}
Notice that tree has separated boundary.

A degree sequence $\pi=(d_1,\dots,d_n)$ of non-negative integers is called a \emph{degree sequence} if there exists a graph $G$ on $n$ vertices for which $d_1,\dots,d_n$ are the degrees of its vertices.
Let $\mathcal{T}_\pi$ denote the set of trees with degree sequence $\pi$.
We show some techniques of rearranging the edges as follows.
\begin{lemma}[Switching {~\cite[Lemma 5]{biyikouglu2007faber}}]\label{lem:Switching1}
    Let $T$ be a tree with boundary in some class $\mathcal{T}_\pi$. 
    Let $v_1, v_2$ be two vertices in $T$ with $v_1 \not\sim v_2$.
    Let $v_1u_1, v_2u_2$ be edges such that $u_2 \in P_T(v_1,v_2)$ and $u_1 \notin P_T(v_1,v_2)$, see~\cref{fig:switching1}. 
    Then by replacing  $v_1u_1$ and $v_2u_2$ by $v_1v_2$ and $u_1u_2$, we obtain a new tree $T'$ which is also contained in $\mathcal{T}_\pi$ with the same set of boundary vertices. Moreover, we find for a function $f \in \RR^{\Omega}$
    \begin{align}\label{R_T_1}
        R_{T'}(f) \leq R_T(f) 
    \end{align}
    whenever $\widehat{f}(v_1) \geq \widehat{f}(u_2)$ and $\widehat{f}(v_2) \geq \widehat{f}(u_1)$. The inequality~\eqref{R_T_1} is strict if both inequalities are strict.
\end{lemma}
\begin{figure}[htbp]
    \centering
      \begin{tikzpicture}[    
x=1.00mm, y=1.00mm,    
inner xsep=0pt, inner ysep=0pt,    
outer xsep=0pt, outer ysep=0pt  
]
\path[line width=0mm] (40,40) rectangle +(120,40);
\definecolor{L}{rgb}{0,0,0}
\definecolor{F}{rgb}{0,0,0}
\path[line width=0.30mm, draw=L, fill=F] (60,60) circle (1.00mm); 
\path[line width=0.30mm, draw=L, fill=F] (80,60) circle (1.00mm); 
\path[line width=0.30mm, draw=L, fill=F] (120,60) circle (1.00mm);
\path[line width=0.30mm, draw=L, fill=F] (140,60) circle (1.00mm);
\path[line width=0.30mm, draw=L] (50,60) -- (60,60);
\path[line width=0.30mm, draw=L] (80,60) -- (120,60);
\path[line width=0.30mm, draw=L] (140,60) -- (150,60);
\path[line width=0.30mm, draw=L, dash pattern=on 1.5mm off 1.5mm] (60,60) -- (80,60);
\path[line width=0.30mm, draw=L, dash pattern=on 1.5mm off 1.5mm] (120,60) -- (140,60);
\path[line width=0.30mm, draw=L] (60,60) to[out=60,in=120] (120,60);
\path[line width=0.30mm, draw=L] (80,60) to[out=-60,in=-120] (140,60);
\node at (60,56){$u_1$};
\node at (80,56){$v_1$};
\node at (120,56){$u_2$};
\node at (140,56){$v_2$};
\end{tikzpicture}
    \caption{Switching :  $v_1u_1$ and $v_2u_2$ are replaced by  $v_1v_2$ and $u_1u_2$.}
    \label{fig:switching1}
\end{figure}
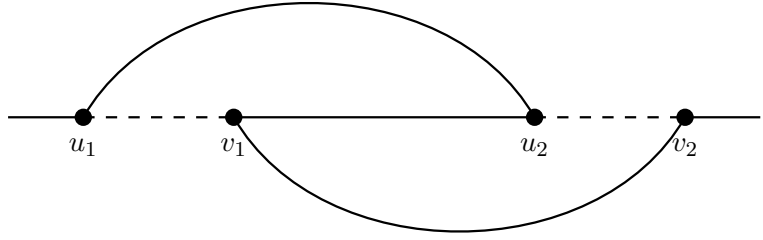
\begin{lemma}[{\cite[Lemma 6]{biyikouglu2007faber}}]\label{lem:strict Switching1}
    Let $T$ be a tree with boundary in some $\mathcal{T}_\pi$ and let $T'$ be a tree obtained from $T$ by applying Switching  as defined in~\cref{lem:Switching1}. 
    If $f$ is a positive eigenfunction of $\lambda_1(T)$, then
    $\lambda_1(T') \leq \lambda_1(T)$ whenever $\widehat{f}(v_1) \ge \widehat{f}(u_2)$ and $\widehat{f}(v_2) \ge \widehat{f}(u_1)$. Moreover, $\lambda_1(T') < \lambda_1(T)$ if at least one of these two inequalities is strict.
\end{lemma}

\begin{figure}[htbp]
    \centering
\begin{tikzpicture}[
    x=1.00mm, y=1.00mm,
    inner xsep=0pt, inner ysep=0pt,
    outer xsep=0pt, outer ysep=0pt
]
\path[line width=0mm] (40,40) rectangle +(120,40);

\definecolor{L}{rgb}{0,0,0}
\definecolor{F}{rgb}{0,0,0}
\path[line width=0.30mm, draw=L, fill=F] (80,60) circle (1.00mm); 
\path[line width=0.30mm, draw=L, fill=F] (130,60) circle (1.00mm);
\path[line width=0.30mm, draw=L, fill=F] (80,45) circle (1.00mm); 

\path[line width=0.30mm, draw=L] (50,60) -- (150,60);

\path[line width=0.30mm, draw=L] (80,35) -- (80,45);

\path[line width=0.30mm, draw=L, dash pattern=on 1.5mm off 1.5mm] (80,60) -- (80,45);

\path[line width=0.30mm, draw=L] (80,45) -- (130,60);

\node at (80,63){$v_1$};
\node at (130,63){$v_2$};
\node at (76,45){$u$};
\end{tikzpicture}
    \caption{Shifting: $uv_1$ is replaced by $uv_2$.}
    \label{fig:shifting}
\end{figure}
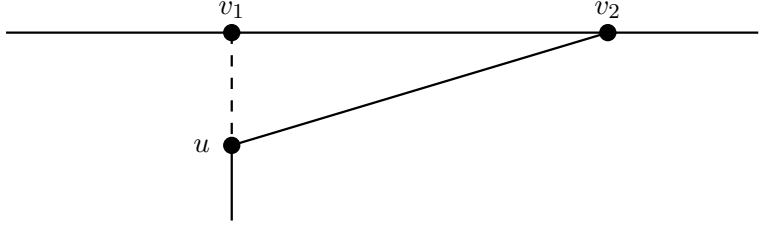

\begin{figure}[htbp]
    \centering
    \begin{tikzpicture}[    
x=1.00mm, y=1.00mm,    
inner xsep=0pt, inner ysep=0pt,    
outer xsep=0pt, outer ysep=0pt  
]
\path[line width=0mm] (40,40) rectangle +(120,40);
  \definecolor{L}{rgb}{0,0,0}  
\definecolor{F}{rgb}{0,0,0}
\path[line width=0.30mm, draw=L, fill=F] (60,60) circle (1.0mm);  
\path[line width=0.30mm, draw=L, fill=F] (100,60) circle (1.0mm); 
\path[line width=0.30mm, draw=L, fill=F] (140,60) circle (1.0mm); 
  \path[line width=0.30mm, draw=L, fill=white] (90,44) circle (1.0mm);  
\path[line width=0.30mm, draw=L, fill=white] (110,44) circle (1.0mm); 
\path[line width=0.30mm, draw=L] (50,60) -- (60,60);  
\path[line width=0.30mm, draw=L, dash pattern=on 1.5mm off 1.5mm] (60,60) -- (100,60);  
\path[line width=0.30mm, draw=L] (100,60) -- (140,60);  
\path[line width=0.30mm, draw=L] (140,60) -- (150,60);
\path[line width=0.30mm, draw=L] (60,60) to[out=60,in=120] (140,60);
\path[line width=0.30mm, draw=L] (100,60) -- (90,45);  
\path[line width=0.30mm, draw=L] (100,60) -- (110,45);
\node at (100,45){$\cdots$};
\node at (60,56){$v_1$};  
\node at (100,56){$u$};  
\node at (140,56){$v_2$};
  \end{tikzpicture}
    \caption{Jumping:  $uv_1$ is replaced by  $uv_2$.}
    \label{fig:jumping}
\end{figure}
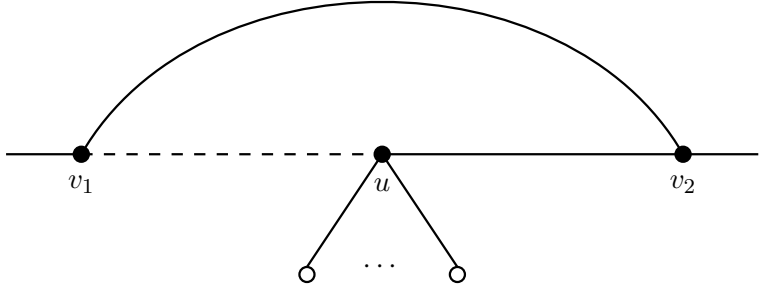

Denote the geodesic path from $u$ to $v$ in the tree $T$ by $P_T(u, v)$.
\begin{lemma}[Shifting{~\cite[Lemma 7]{biyikouglu2007faber}}]\label{lem:Shifting}
    Let $T$ be a tree with boundary in $\mathcal{T}^{(n,k)}$. 
    Let $v_1, v_2$ be two vertices in $T$.
    Let $uv_1$ be an edge with $u \notin P_T(v_1,v_2)$ (see~\cref{fig:shifting}). 
    By replacing $uv_1$ by $uv_2$ we obtain a new tree $T'$ which is also contained in $\mathcal{T}^{(n,k)}$. 
    If $v_2$ is an interior vertex and $d_T(v_1)\ge 3$, then the number of boundary vertices remains unchanged. 
    Moreover, for any positive function $f\in \RR^{\Omega}$ we have
    \begin{align*}
        R_{T'}(f) \leq R_T(f)
    \end{align*}
    if $\widehat{f}(v_1)\geq \widehat{f}(v_2)\geq \widehat{f}(u)$. 
    The inequality is strict if $\widehat{f}(v_1)> \widehat{f}(v_2)$.
\end{lemma}
Define $\mathcal{V}_G \coloneqq \{ v \in \Omega: N_G(v) \cap B \neq \emptyset \}$.
\begin{lemma}[Jumping]\label{lem:Jumping}
    Let $T$ be a tree with boundary in $\mathcal{T}^{(n,k)}$. 
    Let $v_1, v_2$ be two interior vertices in $T$ with $v_1 \not\sim v_2$.
    Let $uv_1$ be an edge with $u\in P_T(v_1,v_2)\cap \mathcal{V}_T$, see~\cref{fig:jumping}.                                      
    By replacing $v_1u$ by $v_1v_2$ we obtain a new tree $T'$ which is also contained in $\mathcal{T}^{(n,k)}$. 
    Then for any positive function $f\in \RR^{\Omega}$ we have
    \begin{align*}
        R_{T'}(f) \leq R_T(f)
    \end{align*}
    if $f(v_1)\geq f(v_2)\geq f(u)$. 
    The inequality is strict if $f(v_2)> f(u)$.
\end{lemma}
\begin{proof}
We need to compute the effects of removing and inserting edges and get
    \begin{align*}
    R_{T'}(f) - R_T(f) = &(f(v_1)-f(v_2))^2  - (f(v_1)-f(u))^2 \\
    = &(f(u)-f(v_2)) (2f(v_1)-f(u)-f(v_2)) \leq 0,
    \end{align*}
where the last inequality is strict if $f(v_2)> f(u)$.
Thus the lemma follows.
\end{proof}

\section{Proofs}

\subsection{\texorpdfstring{Faber–Krahn property in $\mathcal{T}(n,m)$}{Faber–Krahn property in \mathcal{T}(n,m)}}

First, we compare $\lambda_1$ of trees $T(p, 2, b)$ with different lengths and fixed order.
\begin{lemma}\label{D increases-eigenvalue decreases}
    Let $p > \ell \geq 0$ and $b \geq 2$ be integers. 
    Then
    \begin{align*}
        \lambda_1(T(p-\ell, 2, b+\ell)) \geq \lambda_1(T(p, 2, b)).
    \end{align*}
    The equality holds if and only if $\ell = 0$.
\end{lemma}
\begin{proof}
    Let $f$ be an eigenfunction of $\lambda_1(T(p-\ell, 2, b+\ell))$.
    Define a test function on $T(p,2,b)$ by
    \begin{align*}
        g(v)=
        \begin{cases} 
            f(u_i) & v = v_i, 1 \leq i \leq p - \ell+2, \\
            0 & \text{otherwise}.
        \end{cases}
    \end{align*}
    If $\ell \geq 1$,
    then we have
    \begin{align*}
        \lambda_1(T(p,2,b)) &\leq \frac{\sum_{i=1}^{p - \ell + 1} ( g(v_i) - g(v_{i+1}))^2 + g^2(v_{p - \ell+2}) + g^2(v_1)}{\sum_{i=1}^{p - \ell+2} g^2(v_i)} \\
        &
        < \frac{\sum_{i=1}^{p - \ell + 1} ( f(u_i) - f(u_{i+1}))^2 + (b +\ell- 1) f^2(u_{p - \ell+2}) + f^2(u_1)}{\sum_{i=1}^{p - \ell+2} f^2(u_i)} \\
        &= \lambda_1(T(p-\ell,2,b+\ell)).
    \end{align*}
    Therefore, $\lambda_1(T(p-\ell, 2, b+\ell)) = \lambda_1(T(p, 2, b))$ only holds for $\ell=0$.
\end{proof}

Let $P_{\ell}$ be the path $v_0 \sim \cdots \sim v_{\ell-1}$ of length $\ell-1$.
\begin{theorem}[{\cite[Theorem 1.3.]{lin2025estimatesdirichleteigenvaluegraphs}}]\label{thm:lower_bound_2}
    Let $T = (V, E)$ be a tree with leaves as boundary. 
    Let $r$ be the inscribed radius of $T$. 
    Suppose $|V| \geq 3$. 
    Then
    \begin{align*}
        \lambda_1(T) \geq 4 \sin^2 \frac{\pi}{4r+2} \geq \frac{1}{r^2}  
    \end{align*}
    The first equality holds if $T$ is the path graph $P_{2r+2}$.
\end{theorem}

Now we are ready to give the proof of \cref{cor:m_extremal tree}.
\begin{proof}[\bf{{Proof of~\cref{cor:m_extremal tree}}}]
    Let $T$ be a tree with order $n$, interior $\Omega$, diameter $D$ and matching number $m$.
    Clearly, $n \geq 2m$.
    We divide the proof into two cases. 
    \begin{case}\label{thm1:case 1}
     $n \geq 2m+1$.
    
    \begin{subcase}
    $D \leq 2m-2$.
    
    Clearly, $P_{2m} \subseteq T(2m-3, 2, n+1-2m)$.
    By~\cref{thm:subgraph_larger_original graph}, we have
    \begin{align*}
    \lambda_1(T(2m-3, 2, n+1-2m))\leq \lambda_1(P_{2m}) = 2(1-\cos\frac{\pi}{2m-1}).
    \end{align*}
    The first inequality is strict by~\cref{lem:simple eigenvalue-positive}.
    
    If $D \leq 2m-2$, then by~\cref{thm:lower_bound_2} we have
    \begin{align*}
        \lambda_1(T) \geq 4\sin^2\frac{\pi}{4r + 2} \geq 4\sin^2\frac{\pi}{4(m - 1)+2} = 2(1 - \cos\frac{\pi}{2m - 1}).
        \end{align*}
    Therefore, we have $\lambda_1(T) > \lambda_1(T(2m-3, 2, n+1-2m))$.
    \end{subcase}
   
    \begin{subcase}
    $D \geq 2m - 1$.
        
    Take a path $P$ of length $D$.
    If $D \geq 2m+1$, then $P$ has at least $2m+2$ vertices, which implies that the matching number of $P$ is at least $m+1$. 
    This yields a contradiction. 
        
    Next, we consider $D=2m-1$ and $D=2m$.
    
    \begin{claim*}
        For $D=2m-1$ or $D=2m$, $|\Omega| = D-1$.
    \end{claim*}
    
    \begin{proof}
        Take a maximum matching of $P$, which is denoted by $M_P$.
        Clearly, $\card{M_P}=m$.
        Suppose that there exists an interior vertex $v' \notin P$. 
        Then there exists a neighbor  $u'\notin P$ of $v'$. 
        Since $v'$ and $u'$ are uncovered by $M_P$, we can construct a new matching $M = M_P \cup \{v'u'\}$ such that $\card{M}=m+1$, which contradicts the fact that the matching number of $T$ is $m$. 
        Therefore, all the interior vertices are on $P$, which implies that $\card{\Omega}=D-1$.
    \end{proof}
        By~\cref{Klobürštel theorem}, $T(2m-4,2, n+2-2m)$ is the unique tree with the lowest first Dirichlet eigenvalue in $\mathcal{T}(n,m)$ with $|\Omega|=2m-2$, and $T(2m-3,2, n+1-2m)$ is the unique tree with the lowest first Dirichlet eigenvalue in $\mathcal{T}(n,m)$ with $|\Omega|=2m-1$.
        Applying~\cref{D increases-eigenvalue decreases}, we have $\lambda_1(T(2m-4,2, n+2-2m)) > \lambda_1(T(2m-3,2, n+1-2m))$.
        \end{subcase}
    \end{case}
    \begin{case}
    $n = 2m$.
    
    Here $D\leq 2m-1$.
    If $D= 2m-1=n-1$, then $T\cong P_{2m}$.
    If $D\leq 2m-2$, then by~\cref{thm:lower_bound_2}
    \begin{align*}
        \lambda_1(T) \geq 4\sin^2\frac{\pi}{4r + 2} \geq 4\sin^2\frac{\pi}{4(m - 1)+2} = 2(1 - \cos\frac{\pi}{2m - 1})=\lambda_1(P_{2m}).
        \end{align*}
        The last equality holds if and only if $D=2r = 2m-2$, which also implies that $|\Omega|=n-3$.
        By~\cref{Klobürštel theorem}, $T(2m-5,2, 3)$ is the unique tree with the lowest first Dirichlet eigenvalue in $\mathcal{T}(n,m)$ with $|\Omega|=2m-2$.
         Applying~\cref{D increases-eigenvalue decreases}, we have $\lambda_1(T(2m-5,2, 3)) > \lambda_1(T(2m-4,2, 2))$.
    \end{case}
    The proof is complete.
\end{proof}

\subsection{\texorpdfstring{Faber–Krahn property in $\mathcal{T}(n,m, b)$}{Faber–Krahn property in \mathcal{T}(n,m, b)}}

Recall that $\mathcal{V}_G \coloneqq \{ v \in \Omega: N_G(v) \cap B \neq \emptyset \}$.
The following lemma shows that there exists a maximum matching that contains some disjoint pendant edges. 
\begin{lemma}\label{lem:uv-maximum-matching}
    Let $T$ be a tree with matching number $m$. 
    Then for any $|\mathcal{V}_T|$ disjoint pendant edges, there exists a maximum matching of $T$ that contains them.
\end{lemma}
\begin{proof}
    Let $|\mathcal{V}_T|=s$.
    Label the vertices in $ \mathcal{V}_T $ as $ v_1, v_2, \dots, v_s $, and choose an arbitrary vertex $ u_i \in N_T(v_i) \cap B_T $ for each $ v_i $.
    Denote $\{v_1u_1, \cdots, v_su_s \}$ by $S$.
    Take a maximum matching $M_0$ of $T$. 
    If $S\subseteq M_0$, we are done. 
    Otherwise, there exists some $i$ ($1 \leq i \leq s$) such that $u_i v_i \notin M_0$.
    If  $M_0$ does not cover $v_i$, then
    $\lvert M_0\cup\{u_iv_i\}\rvert = m+1$, a contradiction. 
    Hence, the vertex $v_i$ is covered by the edge $v_iw_i\in M_0$, where $w_i\neq u_i$.
    Define
    \begin{align*}
        M_1 \coloneqq (M_0\setminus\{v_iw_i\})\cup\{u_iv_i\}.
    \end{align*}
    It is clear that $M_1$ is also a maximum matching, and it contains $u_iv_i$.
    Since $S$ is finite, we can repeat the operation above and obtain a sequence $M_0,M_1,\cdots, M_c$ for some finite $c$ such that $S \subseteq M_c$.
    It is clear that $\card{M_c}=\card{M_0}$.
    Therefore, $M_c$ is the desired matching.
\end{proof}

Next, we give a lower bound for $\card{\mathcal{V}_T}$.
\begin{lemma}\label{lem:s_t}
Let $T \in \mathcal{T}(n,m,b)$. 
Then $\card{\mathcal{V}_T} \geq 2m+b-n$.
\end{lemma}
\begin{proof}
    Let $\card{\mathcal{V}_T} =s$ and $t=2m+b-n$.
    By~\cref{lem:uv-maximum-matching}, there exists a maximum matching $M$  that covers $s$ leaves. 
    Since any two edges in $M$ are vertex-disjoint, the remaining $b-s$ leaves cannot be covered by $M$, which means $b-t = n - 2m \geq b-s$. 
    It follows that $s \geq t$.
\end{proof}

Let $t=2m+b-n$.
The following two lemmas give the range of $t$.
\begin{lemma}\label{lem:tree-n-bound}
    Let $T$ be a tree with $n\geq 2$ vertices, $b$ leaves, and matching number $m$. 
    Then $n\leq 2m+b-1$.
\end{lemma}
\begin{proof}
    We proceed by induction on $n$.
    For $n=2,3$, the inequality holds clearly.
    Assume that $n\geq 4$ and the result holds for all trees with the number of vertices fewer than $n$. 
    Let $u$ be a leaf and $v$ be a neighbor  of $u$. 
    Define $T' \coloneqq T \setminus \{u,v\}$.
    Denote the order, the number of leaves, and matching number of $T'$ by $n',b',m'$ respectively.
    Set
    \begin{align*}
        x& =\card{\{ w\in N_T(v)\setminus\{u\}: d_w=1\}},\\
        y& =\card{\{ w\in N_T(v): d_w=2\}},\\
        z& =\card{\{ w\in N_T(v): d_w\geq 3\}}.
    \end{align*}
    Then $n' = n-2$ and $b' = b - 1 - x + y$.
    \begin{claim*}
        $m' = m-1$
    \end{claim*} 
    \begin{proof}
        By~\cref{lem:uv-maximum-matching}, there exists a maximum matching $M$ in $T$ which contains $uv$.
        Denote $M_1 \coloneqq M\setminus\{uv\}$. 
        Then $M_1$ is a matching of $T'$,  which implies $m' \geq |M_1| = |M|-1=m-1$.
        Conversely, a matching of $T$ can be obtained by $M' \cup \{uv\}$ where $M'$ is a maximum matching of $T'$, implying that $m \geq 1 + |M'|=1+m'$.
    \end{proof}
    Applying the inductive hypothesis to each nontrivial component $T_i$ of $T'$, we have $n_i \leq 2m_i + b_i - 1$.
    Summing over all $y+z$ components gives
    \begin{align*}
        \sum_{i=1}^{y+z} n_i \leq \sum_{i=1}^{y+z} ( 2m_i+b_i - 1) 
        = \sum_{i=1}^{y+z} (2m_i+b_i) - (y+z).
    \end{align*} 
    Note that $n'= \sum_{i=1}^{y+z} n_i + x$, $m'=\sum_{i=1}^{y+z} m_i$, and $b'=\sum_{i=1}^{y+z} b_i$. 
    Combining $n' = n-2$, $b' = b - 1 - x + y$ and $m' = m-1$,
    we have
    \begin{align*}
         n &= n'+2 =\sum_{i=1}^{y+z} n_i + x+2 \leq \sum_{i=1}^{y+z} (2m_i+b_i) - (y+z) + x +2\\
         &= 2m'+ b' - (y+z) + x+2=2m'+(b - 1 - x + y) - (y+z) + x+2\\
         &= 2m'+ b - z +1 \leq 2m+b-1.   \qedhere
    \end{align*} 
\end{proof}

The lemma above implies that $t\geq 1$.
\begin{lemma}\label{lem:t_upper bound}
    Let $n\geq3, b\geq2, m\geq1$ be integers.
    Let $t=2m+b-n$.
    If $T \in \mathcal{T}(n,m,b)$, then $t \leq \min\{b,m\}$.
\end{lemma}
\begin{proof}
    Let $M$ be a maximum matching of $T$.
    The number of vertices uncovered by $M$ is $n - 2m$, which implies that $n - 2m = b - t \geq 0$.
    Define a mapping $\phi: M \to \Omega$ as follows. 
    For each edge $e \in M$, map it to one of its endpoints which is an interior vertex.
    This is feasible as every edge in a tree of order at least three has at least one interior endpoint.
    Therefore, $\phi$ is injective, which implies that $|M| = m \leq |\Omega| = n - b = 2m - t$.
\end{proof}

Before our proof of~\cref{thm:m_b_extremal tree}, we recall two facts as follows.
\begin{theorem}[{\cite[Corollary 1.5.]{lin2025estimatesdirichleteigenvaluegraphs}}]\label{tree-up-bound}
    Let $T$ be a tree with the set of leaves as boundary $B$. 
    Then
    \begin{align*}
        \lambda_{1}(T)\leq \frac{|B|}{|\Omega|}.
    \end{align*}   
    The equality holds if and only if $|B|$ is divisible by $|\Omega|$ and each interior vertex is adjacent to $|B|/|\Omega|$ leaves.
\end{theorem}

\begin{lemma}[{\cite[Theorem 4.2.2]{Horn_Johnson_2012}}]\label{lem:extreme_eigenfunction}
    Let $M$ be a $n \times n$ real symmetric matrix. 
    Let $\lambda_1 < \lambda_2 < \ldots < \lambda_k$ be the distinct eigenvalues of $M$. 
    If $x^\top M x = \lambda_1 x^\top x$ for $x \in \RR^n$, then $Mx = \lambda_1 x$. 
    If $x^\top M x = \lambda_k x^\top x$ for $x \in \RR^n$, then $Mx = \lambda_k x$. 
\end{lemma}

We now proceed to prove our main theorem.

\begin{proof}[\bf{{Proof of~\cref{thm:m_b_extremal tree}}}]

By~\cref{lem:tree-n-bound} and~\cref{lem:t_upper bound}, we have
$1\leq t \leq \min\{b,m\}$.
Choose a tree $T \in \mathcal{T}(n,m,b)$.
By~\cref{lem:s_t}, we have $\card{\mathcal{V}_T} \geq t$.
Denote the number of interior vertices of $T$ by $k=n-b$.
\setcounter{case}{0} 

\begin{case}\label{case1:m_b_extremal tree}
$m >t \geq 2$
    
It follows that $k \geq 3$.
Let $f\in\mathbb{R}^\Omega$ be a positive eigenfunction of the first Dirichlet eigenvalue of $T$ and let $\widehat{f}$ be an extension of $f$ to the whole tree by assigning 0 to the leaves.
For convenience, we still use $f$ to denote $\widehat{f}$.
Assume that the vertices of $T$, $V = \{v_1, v_2, \ldots, v_{k}, \ldots, v_n\}$, are numbered such that $f(v_i) \geq f(v_j)$ if $i \leq j$. 
        
Through a series of operations, which would not change $\Omega$ and $B$, we can construct the desired tree $T^*=T(2m-2t,t,b)$ with degree sequence $\pi^*$ as follows.
We start with the vertex $v_1$.
Denote the neighbors of $v_1$ by $v_{1,1}, v_{1,2}, \cdots, v_{1,d'_T(v_1)}, \cdots, v_{1,d_T(v_1)}$, where 
$v_{1,1}, v_{1,2}, \cdots, v_{1,d'_T(v_1)} \in N_T(v_1) \cap \Omega$ and $d'_T(v_1) = \card{N_T(v_1) \cap \Omega}$.
If $v_1 \sim v_2$, there is nothing to do, and let $T_1$ denote $T$. 
Now, consider $v_1 \not\sim v_2$.
If there exists a vertex $v_1' \in N_T(v_1)$ such that $f(v_1') = f(v_2)$, 
then we exchange the positions of $v_1'$ and $v_2$ in the ordering of $V$ (and update the indices of the vertices) and let $T_1$ denote $T$. 
If we have $f(v) < f(v_2)$ for every $v \in N_T(v_1)$, 
then there exists a neighbor  $v_{1,1}$ of $v_1$ such that $v_{1,1} \in P_T(v_1,v_2)$ and $f(v_{1,1}) < f(v_2)$.
There are two cases for $v_{1,1}$.
\begin{itemize}
    \item $v_{1,1} \in \mathcal{V}_T$.
    
Using Jumping step, we can replace $v_1v_{1,1}$ by $v_1 v_2$, and obtain a tree $T_1$ with $\card{\mathcal{V}_{T_1}} \geq t$.
By~\cref{lem:Jumping}, we have $R_{T_1}(f)< R_{T}(f)$.

\item $v_{1,1} \notin \mathcal{V}_T$.

Using Switching, we can exchange $v_{1,1}$ with $v_2$, and then $v_2$ becomes a neighbor of $v_1$.
Since $v_2 \notin B$, there exists a neighbor  $v_{2,1}$ of $v_2$ such that
$v_{2,1} \notin P_T(v_1,v_2)$.
By the construction above, we have $f(v_1)\geq f(v_{2,1})$ and $f(v_2) > f(v_{1,1})$.
We replace edges $v_1v_{1,1}$ and $v_2v_{2,1}$ by $v_1v_2$ and $v_{1,1}v_{2,1}$, and obtain a tree $T_1$ with $\card{\mathcal{V}_{T_1}} \geq t$.
By~\cref{lem:Switching1}, we have $R_{T_1}(f)\leq R_{T}(f)$.
\end{itemize}

Next, we obtain another new tree $T_2$ such that $v_1 \sim v_3$. 
If $v_1 \sim v_3$ already holds, then there is nothing to do, and let $T_2$ denote $T_1$. 
If $v_1 \not\sim v_3$, then there exists a vertex $v_{1,2} \in P_{T_1}(v_1,v_3) \cap N_{T_1}(v_1)$.  
We consider two cases for $v_{1, 2}$.

First, if $f(v_{1,2}) \leq f(v_3)$, then by the same procedure as above, we can exchange another vertex $v_{1,2}$ adjacent to $v_1$ (if necessary) with the corresponding vertex $v_3$, and obtain a tree $T_2$ with $\card{\mathcal{V}_{T_2}} \geq t$. 
By~\cref{lem:Switching1} and~\cref{lem:Jumping}, we have $R_{T_2}(f)\leq R_{T_1}(f)$.

Second, if $f(v_{1,2}) > f(v_3)$, then $v_{1,2}=v_2$.
If there exists $v_{1,2}' \neq v_2 \in N_{T_1}(v_1)\setminus B$, then we use Switching.
Since $T_1$ is connected, there exists a neighbor  $v_{3,1}$ of $v_3$ such that $v_{3,1} \in P_{T_1}(v_1,v_3)$.
By the construction, we have $f(v_1)\geq f(v_{3,1})$ and $f(v_3) > f(v_{1,2}')$.
We replace $v_1v_{1,2}' $ and $v_3v_{3,1}$ by $v_1v_3$ and $v_{1,2}' v_{3,1}$, and obtain a new tree $T_2$ with $\card{\mathcal{V}_{T_2}} \geq t$.
By~\cref{lem:Switching1}, we have $R_{T_2}(f)\leq R_{T_1}(f)$. 
If such a vertex $v_{1,2}'$ does not exist, which implies that $v_2$ is the only interior neighbor of $v_1$, then we divide the discussion into following subcases based on $v_{3,1}$.
Select a vertex $v'\in N_{T_1}(v_1)\cap B$. 
\begin{itemize}
    \item $v_{3,1} \notin \mathcal{V}_{T_1}$.
    
    Clearly, $f(v_1)\geq f(v_{3,1})$ and $f(v_3) > f(v')$.
    We replace $v_1v' $ and $v_3v_{3,1}$ by $v_1v_3$ and $v' v_{3,1}$, and obtain a new tree $T_2$ with $\card{\mathcal{V}_{T_2}} \geq t$.
    By~\cref{lem:Switching1}, we have $R_{T_2}(f)\leq R_{T_1}(f)$.
    
    \item  $v_{3,1} \in \mathcal{V}_{T_1}$ and $f(v_{3,1}) \geq f(v_{k-t+1})$.
    
    By the same procedure as the case $v_{3,1} \notin \mathcal{V}_{T_1}$, 
    we obtain $T_2'$ such that $v_{3,1} \sim v'$. 
    For $v_{k-t+1}, v_{k-t+2}, \cdots, v_{k}$, if there exists $v_i$  for $k-t+1 \leq i \leq k$ such that $v_i \notin \mathcal{V}_{T_2'}$, then we use Shifting to replace $v_{3,1}v'$ by $v_iv'$.
    If $v_i \in \mathcal{V}_{T_2'}$ for any $k-t+1 \leq i \leq k$, then we use Shifting to replace $v_{3,1}v'$ by $v_k v'$.
    We obtain a new tree $T_2$ with $\card{\mathcal{V}_{T_2}} \geq t$.
    By~\cref{lem:Switching1} and~\cref{lem:Shifting}, we have $R_{T_2}(f)\leq R_{T_2'}(f)$.
    
    \item $v_{3,1} \in \mathcal{V}_{T_1}$ and any $v$ with $f(v_{k-t+1}) \geq f(v) > f(v_{3,1})$ satisfies that  $v \in \mathcal{V}_{T_1}$.
    
    By the same procedure as $v_{3,1} \notin \mathcal{V}_{T_1}$, 
    we also obtain $T_2'$ such that $v_{3,1} \sim v'$. 
    Denote $v_{3,1}$ by $v_j$, where $k-t+1< j \leq k$.
    According to the assumption, $v_{k-t+1}, \cdots, v_j$ are adjacent to some leaves.
    If there exists $v_i$ for $j\leq i\leq k$ such that $v_i \notin \mathcal{V}_{T_2'}$, then we use Shifting to replace $v_{3,1}v'$ by $v_iv'$.
    If $v_i \in \mathcal{V}_{T_2'}$ for any $j \leq i \leq k$, then we use Shifting to replace $v_{3,1}v'$ by $v_k v'$.
    We obtain a new tree $T_2$ with $\card{\mathcal{V}_{T_2}} \geq t$.
    By~\cref{lem:Switching1} and~\cref{lem:Shifting}, we have $R_{T_2}(f)\leq R_{T_2'}(f)$.

    \item $v_{3,1} \in \mathcal{V}_{T_1}$ and there exists some vertex $v$ with $f(v_{k-t+1}) \geq f(v) > f(v_{3,1})$ that satisfies $v \notin \mathcal{V}_{T_1}$.
    
    We replace $v_3v_{3,1}$ by $v_3v$, and obtain $H$ with $\card{\mathcal{V}_H} \geq t$.
    Note that $H$ can be disconnected when $v_3 \in P_{T_1}(v_{3,1},v)$.
    We do not consider the new leaf as boundary vertex.
    Since $k-t+1=2m-2t+1 \geq 2m-2(m-1)+1=3$, it follows that $f(v_3) \geq f(v) > f(v_{3,1})$.
    Thus, we have
    \begin{align*}
        R_{H}(f) - R_{T_1}(f) =&(f(v_3)-f(v))^2-(f(v_3)-f(v_{3,1}))^2\\
        =&(f(v_{3,1})-f(v))(2f(v_3)-f(v_{3,1})-f(v))] < 0.
    \end{align*}
   
    If $v_3 \notin  P_{T_1}(v_{3,1},v)$, then $H$ is still a tree.
    Select a vertex $v'\in N_H(v_1)\cap B_H$, and we have $f(v_1)\geq f(v)$ and $f(v_3) > f(v')$.
    Apply Switching  by replacing $v_1v'$ and $v_3v$ by $v_1v_3$ and $v'v$, and obtain a new tree $T_2$ with $\card{\mathcal{V}_{T_2}} \geq t$ such that $R_{T_2}(f)\leq R_{H}(f)$ by~\cref{lem:Switching1}.
   
    If $v_3 \in  P_{T_1}(v_{3,1},v)$, then there exists a multiple edge or a cycle in $H$. 
    We replace $v_3v$ and $v_1v'$ by $v_3v_1$ and $vv'$ to obtain $T_2$ with $\card{\mathcal{V}_{T_2}} \geq t$. 
    This operation ensures that $T_2$ is a tree since $v_2$ is the only interior neighbor of $v_1$, which implies that $v_1$ is not in the cycle of $H$.
    Since $f(v_1) \geq f(v)$ and $f(v_3) > f(v')$, we have
    \begin{align*}
        R_{T_2}(f) - R_{H}(f) = &[(f(v_3)-f(v_1))^2 + (f(v)-f(v'))^2] \\ & - [(f(v_1)-f(v'))^2 + (f(v_3)-f(v))^2] \\
    = &2(f(v')-f(v_3))\cdot (f(v_1)-f(v)) \leq 0.
    \end{align*}
    \end{itemize}
    
Now we obtain a tree $T_2$ such that $v_3\sim v_1$ and $R_{T_2}(f)\leq R_{T_1}(f)$. 
If $d_{T_2}(v_1)=2$, there is nothing to do, and let $T^{(1)}$ denote $T_2$.
Otherwise, let $\{v_2,v_3\}=\{v_{1,1},v_{1,2}\}$, and we can use Shifting to replace $v_1v_{1,3}$, $v_1v_{1, 4}$, $\cdots$, $v_1v_{1,d'_T(v_1)}$ with $v_2v_{1,3}$, $v_2v_{1, 4}$, $\cdots$, $v_2v_{1,d'_T(v_1)}$.
For $v_{1,d'_T(v_1)+1}, \cdots, v_{1,d_T(v_1)}$, we reattach them to the interior vertices among $v_{k-t+1}, \cdots, v_{k-1}, v_{k}$. 
For $v_i$ $(k-t+1\leq i\leq k)$, if $v_i$ is not adjacent to any leaf, we shift one of the remaining leaves in $\{v_{1,d'_T(v_1)+1}, \cdots, v_{1,d_T(v_1)}\}$ to it.
If there are still some leaves left after processing all these interior vertices, we shift all the remaining leaves to $v_k$.
After the operations above, we obtain a tree $T^{(1)}$ with $\card{\mathcal{V}_{T^{(1)}}} \geq t$.
As $f(v_1) \geq f(v_2) \geq f(v_{1,3}), \cdots$, $f(v_{1,d'_T(v_1}))$ and $f(v_1) \geq f(v_{k-t+1}),  \cdots, f(v_{k-1}) > 0$, by~\cref{lem:Shifting}, we have 
\begin{align}\label{ali:T_Faber-Krahn property}
    R_{T^{(1)}}(f) \leq R_{T_2}(f).
\end{align}

Assume that we have obtained a tree $T^{(i-1)}$ such that 
$\card{\mathcal{V}_{T^{(i-1)}}} \geq t$, $d_{T^{(i-1)}}(v_{i-1})=2$, $v_1 \sim v_2$, and $v_j \sim v_{j+2}$ for all $1 \leq j \leq i-1$ for some $i \geq 2$. 
In the step $T_1 \to T_2$, we do not require the condition $f(v_1) \geq f(v_2)$. 
Therefore, by the similar argument, we can obtain a tree $T_2^{(i)}$ from $T^{(i-1)}$ such that $R_{T_2^{(i)}}(f) \leq R_{T^{(i-1)}}(f)$, $\card{\mathcal{V}_{T_2^{(i)}}} \geq t$, $v_1 \sim v_2$, and $v_j \sim v_{j+2}$ for all $1 \leq j \leq i$.
By applying the argument similar to the step $T_2 \to T^{(1)}$, 
we obtain the tree $T^{(i)}$ from $T_2^{(i)}$ such that $R_{T^{(i)}}(f) \leq R_{T_2^{(i)}}(f)$, $d_{T^{(i)}}(v_{i}) = 2$, $\card{\mathcal{V}_{T^{(i)}}} \geq t$, $v_1 \sim v_2$, and $v_j \sim v_{j+2}$ for all $1 \leq j \leq i$.

Through the operations above, we obtain a sequence of trees $$T^{(1)}\to T^{(2)} \to T^{(3)} \to \cdots \to T^{(k-t)}$$ such that  $d_{T^{(k-t)}}(v_i) = 2$ for any $1 \leq i \leq k-t$, $\card{\mathcal{V}_{T^{(k-t)}}}= t$, and
\begin{align*}
    R_{T^{(k-t)}}(f) \leq \cdots \leq R_{T^{(2)}}(f) \leq R_{T^{(1)}}(f).
\end{align*}
Moreover, $v_i \in \mathcal{V}_{T^{(k-t)}}$ for $k-t+1 \leq i \leq k$ and $v_{k-t+2} \sim v_{k-t} \sim \cdots  \sim v_{k-t-1} \sim v_{k-t+1}$.
Now, we perform operations on $T^{(k-t)}$.  
If $N_{T^{(k-t)}}(v_{k-t+1})\setminus B = v_{k-t-1}$ and $\card{N_{T^{(k-t)}}(v_{k-t+1})\cap B} = 1$, then there is nothing to do, and let $T^{(k-t+1)}$ denote $T^{(k-t)}$.  
Otherwise $v_{k-t-1}$ must be adjacent to some other vertices, denoted by $v_{k-t+1,1}, v_{k-t+1,2}, \cdots$, and then we can shift these vertices to $v_{k-t+2}$, and obtain a new tree $T^{(k-t+1)}$ with $\card{\mathcal{V}_{T^{(k-t+1)}}} \geq t$.  
Since $f(v_{k-t+1}) \geq f(v_{k-t+2}) \geq f(v_{k-t+1,1}), f(v_{k-t+1,2}), \cdots$, by~\cref{lem:Shifting}, we have
\begin{align*}
    R_{T^{(k-t+1)}}(f) \leq R_{T^{(k-t)}}(f).
\end{align*}
We next perform operations on $T^{(k-t+1)}$ based on three cases.
\begin{itemize}
\item  $v_{k-t+2} \sim v_{k-t+3}$.

There is nothing to do, and let $T_2^{(k-t+2)}$ denote $T^{(k-t+1)}$.
\item $v_{k-t+2} \not\sim v_{k-t+3}$ and there exists a vertex $v_{k-t+2}'\neq v_{k-t} \in N_{T^{(k-t+1)}}(v_{k-t+2})$ such that 
$f(v_{k-t+2}') = f(v_{k-t+3})$.
We exchange the positions of 
$v_{k-t+2}'$ and $v_{k-t+3}$ in the ordering of $V$ (and update the indices 
of the vertices), and let $T_2^{(k-t+2)}$ denote $T^{(k-t+1)}$.
\item $v_{k-t+2} \not\sim v_{k-t+3}$ and $f(v) < f(v_{k-t+3})$ for any $v\neq v_{k-t} \in N_{T^{(k-t+1)}}(v_{k-t+2})$. 

There exists a vertex $v_{k-t+2}' \in N_{T^{(k-t+1)}}(v_{k-t+2}) \cap P_{T^{(k-t+1)}}(v_{k-t+2},v_{k-t+3})$ such that $f(v_{k-t+2}') < f(v_{k-t+3})$.
We replace
$v_{k-t+2}v_{k-t+2}'$ by $v_{k-t+2}v_{k-t+3}$, and obtain a tree 
$T_2^{(k-t+2)}$ with $\card{\mathcal{V}_{T_2^{(k-t+2)}}} \geq t$. 
By~\cref{lem:Jumping}, we have
\begin{align*}
    R_{T_2^{(k-t+2)}}(f) < R_{T^{(k-t+1)}}(f).
\end{align*}
\end{itemize}
If $N_{T_2^{(k-t+2)}}(v_{k-t+2}) \setminus (B \cup \{v_{k-t}, v_{k-t+3}\})$ is nonempty or $\card{N_{T_2^{(k-t+2)}}(v_{k-t+2}) \cap B} \geq 2$, then we apply Shifting to other neighbors of $v_{k-t+2}$ in the same manner as we have done for $v_{k-t+1}$, and obtain a new tree $T^{(k-t+2)}$  with $\card{\mathcal{V}_{T^{(k-t+2)}}} \geq t$ such that
\begin{align*}
    R_{T^{(k-t+2)}}(f) \leq R_{T_2^{(k-t+2)}}(f).
\end{align*}
Otherwise, we directly denote $T_2^{(k-t+2)}$ by $T^{(k-t+2)}$.
With the similar argument, we can use a series of Jumping or Shifting operations to obtain a sequence of trees 
$$T^{(k-t)}\to T^{(k-t+1)}\to T^{(k-t+2)}\to T^{(k-t+3)} \to T^{(k-t+4)} \to \cdots \to T^{(k-1)}$$
 with $\card{\mathcal{V}_{T^{(k-1)}}} = t$ ($T^{(k-1)}$ is shown as~\cref{fig:schematic graph}) such that
\begin{align*}
    \lambda_1(T^*) \leq R_{T^{(k)}}(f) \leq \cdots  \leq R_{T^{(k-t+2)}}(f) \leq R_{T^{(k-t+1)}}(f) \leq R_{T^{(k-t)}}(f). 
\end{align*}
\begin{figure}[htbp]
    \centering
    \begin{tikzpicture}[    
x=1.00mm, y=1.00mm,    
inner xsep=0pt, inner ysep=0pt,    
outer xsep=0pt, outer ysep=0pt,    
node circle/.style={line width=0.30mm, draw=L, fill=F, circle},    
node circle small/.style={line width=0.15mm, draw=L, fill=F, circle}
]
\path[line width=0mm] (41.90,36.47) rectangle +(66.10,26.59);
\definecolor{L}{rgb}{0,0,0}
\definecolor{F}{rgb}{0,0,0}
\path[line width=0.30mm, draw=L] (40,60) circle (1.00mm);
\path[line width=0.30mm, draw=L, fill=F] (40,50) circle (1.00mm);
\path[line width=0.30mm, draw=L, fill=F] (52,50) circle (1.00mm);
\node[inner sep=0pt] at (60,50)  
{\tikz{\fill (0,0) circle (0.5mm);\fill (2mm,0) circle (0.5mm);\fill (4mm,0) circle (0.5mm);}};
\path[line width=0.30mm, draw=L, fill=F] (72,50) circle (1.00mm);
\path[line width=0.30mm, draw=L, fill=F] (84,50) circle (1.00mm);
\path[line width=0.30mm, draw=L, fill=F] (96,50) circle (1.00mm);
\node[inner sep=0pt] at (106,50)  
{\tikz{\fill (0,0) circle (0.5mm);\fill (2mm,0) circle (0.5mm);\fill (4mm,0) circle (0.5mm);}};
\path[line width=0.30mm, draw=L, fill=F] (116,50) circle (1.00mm); 
\path[line width=0.30mm, draw=L, fill=F](128,50) circle (1.00mm);
\path[line width=0.30mm, draw=L, fill=F] (140,50) circle (1.00mm);
\path[line width=0.30mm, draw=L] (128,60) circle (1.00mm);
\path[line width=0.30mm, draw=L] (140,60) circle (1.00mm);
\node[inner sep=0pt] at (150,50)  
{\tikz{\fill (0,0) circle (0.5mm);\fill (2mm,0) circle (0.5mm);\fill (4mm,0) circle (0.5mm);}};
\path[line width=0.30mm, draw=L, fill=F] (160,50) circle (1.00mm);
\path[line width=0.30mm, draw=L, fill=F] (172,50) circle (1.00mm);
\path[line width=0.30mm, draw=L] (168,60) circle (1.00mm);
\path[line width=0.30mm, draw=L] (176,60) circle (1.00mm);
\node at (172,60) {$\cdots$};
\path[line width=0.30mm, draw=L] (160,60) circle (1.00mm);
\node at (40,47){$v_{k-t+1}$};
\node at (52,47){$v_{k-t-1}$};
\node at (72,47){$v_3$};
\node at (84,47){$v_1$};
\node at (96,47){$v_2$};
\node at (116,47){$v_{k-t}$};
\node at (128,47){$v_{k-t+2}$};
\node at (140,47){$v_{k-t+3}$};
\node at (160,47){$v_{k-1}$};
\node at (172,47){$v_{k}$};
\path[line width=0.30mm, draw=L] (40,50) -- (40,59);
\path[line width=0.30mm, draw=L] (40,50) -- (55,50);
\path[line width=0.30mm, draw=L] (69,50) -- (99,50);
\path[line width=0.30mm, draw=L] (113,50) -- (143,50);
\path[line width=0.30mm, draw=L] (128,50) -- (128,59);
\path[line width=0.30mm, draw=L] (140,50) -- (140,59);
\path[line width=0.30mm, draw=L] (157,50) -- (172,50);
\path[line width=0.30mm, draw=L] (160,50) -- (160,59);
\path[line width=0.30mm, draw=L] (172,50) -- (168,59);
\path[line width=0.30mm, draw=L] (172,50) -- (176,59);
\end{tikzpicture}
    \caption{$T^{(k-1)}$}
    \label{fig:schematic graph}
\end{figure}
Actually, $T^*\cong T^{(k-1)}$.
Finally, we show that $T^*$ is the unique tree that has Faber-Krahn property in $\mathcal{T}(n,m,b)$. 
Assume that there exists another tree $T' \ncong T^*$ that has Faber-Krahn property in $\mathcal{T}(n,m,b)$.
Denote $T'$ by $T^{(0)}$.
We do all the operations to $T'$ as above, and get a sequence of trees such that
\begin{align*}
     \lambda_1(T^*) = R_{T^{(k-1)}}(f) \leq \cdots \leq R_{T^{(k-t)}}(f) \leq \cdots \leq R_{T^{(1)}}(f)\leq \lambda_1(T^{(0)}).
\end{align*}
Notice that all the equalities here must hold, and $f$ should be the eigenfunction corresponding to the first Dirichlet eigenvalue for every tree in this sequence.

To lead to a contradiction, we show that it is impossible to perform any operation from $T^{(i-1)}$ to $T^{(i)}$, for $1 \leq i \leq k-1$.
Let $v_{i,1}$ be a neighbor of $v_i$ such that $v_{i,1} \in P_{T^{(i-1)}}(v_i, v_{i+2})$.

First we prove that Switching cannot be applied.
For $T' \to T'_1$, the existence of Switching implies that 
$v_{1,1} \neq v_2$ and $f(v_{1,1}) < f(v_2)$.
Then $\lambda_1(T'_1) <\lambda_1(T')$ by~\cref{lem:strict Switching1}, a contradiction to the Faber-Krahn property of $T'$.
For $k-t+1 \leq i \leq k-1$, we do not need any Switching operation in our proof.
So we only consider $1 \leq i \leq k-t$.
The existence of Switching implies that 
$v_{i,1} \neq v_{i+2}$ and $f(v_{i,1}) \neq f(v_{i+2})$.
If $f(v_{i+2}) > f(v_{i,1})$, then $\lambda_1(T^{(i)}) <\lambda_1(T^{(i-1)})$ by~\cref{lem:strict Switching1}, a contradiction.
If $f(v_{i+2}) < f(v_{i,1})$. 
there exist a neighbor $v_{i+2}'$ of $v_{i+2}$
with $v_{i+2}' \in P_{T^{(i-1)}}(v_i, v_{i+2})$ and a neighbor $v_i'$ of $v_i$ with
$v_i' \notin P_{T^{(i-1)}}(v_i, v_{i+2})$ such that
$f(v_{i+2}) > f(v_i')$ and $f(v_i) \ge f(v_{i+2}')$.
Then $\lambda_1(T^{(i)}) <\lambda_1(T^{(i-1)})$ by~\cref{lem:strict Switching1}, a contradiction.

Next, we shall prove that the Jumping operation does not exist.
Suppose that Jumping exists for $1 \leq i \leq k-t$, which implies that $v_{i,1} \in \mathcal{V}_{T^{(i-1)}}$ and
$f(v_{i+2}) > f(v_{i,1})$.
Then $R_{T^{(i)}}(f) < R_{T^{(i-1)}}(f)$ by~\cref{lem:Jumping}, a contradiction.
Suppose that jumping exists for $k-t+2 \leq i \leq k-2$. 
The proof is analogous to the case $1 \leq i \leq k-t$.

Finally, we show that the Shifting operation cannot be applied.
We state the following claim about the eigenfuntion of $T^*$.
\begin{claim}\label{cla:f_decreases}
    $f(v_j) > f(v_{j+2})$ for $1 \leq j \leq k-t+1$
    and $f(v_j) > f(v_{j+1})$ for $k-t+2 \leq j \leq k-1$.
\end{claim}
\begin{proof}
Clearly, $P_5 \subseteq T^*$.
By~\cref{thm:subgraph_larger_original graph}, we have
\begin{align*}
    \lambda_1(T^*)\leq \lambda_1(P_5) = 2(1 - \cos\frac{\pi}{4}) ~(\cite[\text{Example} \ 2.8]{lin2025estimatesdirichleteigenvaluegraphs})< 1 \leq d_{T^*}(v_k)-1.
\end{align*}
Therefore, we have 
\begin{align*}
    \lambda_1(T^*) f(v_k) = d_{T^*}(v_k)f(v_k) - f(v_{k-1})< (d_{T^*}(v_k)-1)f(v_k),
\end{align*}
which indicates that $f(v_k)<f(v_{k-1})$.

Suppose that $f(v_j) = f(v_{j-1})$ for some $j$ such that $k-t+3 \leq j \leq k-1$.
It follows that
\begin{align*}
    \lambda_1(T^*) f(v_j) = 3f(v_j) - f(v_{j-1})-f(v_{j+1}) < f(v_j),
\end{align*}
which leads to $f(v_j) < f(v_{j+1})$, a contradiction. 
Therefore, we have $f(v_{j+1}) < f(v_{j})$ for any $k-t+2 \leq j \leq k-1$.
Suppose that $f(v_j) = f(v_{j+2})$ for $2 \leq j \leq k-t$.
Then
\begin{align*}
    \lambda_1(T^*) f(v_j) = 2f(v_j) - f(v_{j-2}) - f(v_{j+2}) \leq 0,
\end{align*}
which contradicts $f(v_j)>0$ and $\lambda_1(T^*)>0$.
Therefore, we have $f(v_{j+2}) < f(v_{j})$ for any $2 \leq j \leq k-t$.
For $v_1$, if $f(v_1)=f(v_2)=f(v_3)$, then $\lambda_1(T^*)f(v_1)= 2f(v_1) - f(v_2) -f(v_3) = 0$, again a contradiction.
For $v_{k-t+1}$, we have $f(v_{k-t+1}) \geq f(v_{k-t+2})>f(v_{k-t+3})$.
Therefore, we have $f(v_{j+2}) < f(v_j)$ for $1 \leq j \leq k-t+1$.
\end{proof}
Now we prove that Shifting does not exist.
We only prove for $1 \leq i \leq k-t$ here, and the proof for $k-t+1 \leq i \leq k-1$ is analogous.
Suppose that Shifting exists for $1 \leq i \leq k-t$. 
Then $d_{T^{(i)}}(v_i) \geq 3$. 
Let $v_{-1} = v_2$ and $v_{0} = v_{1}$.
Since $f$ is an eigenfunction of $T^*$ and $T^{(i)}$ and  $\lambda_1(T^*)f(v_i) = \lambda_1(T^{(i)})f(v_i)$, we have
\begin{align*}
    &d_{T^*}(v_i)f(v_i)-f(v_{i-2})- f(v_{i+2})\\
    &= d_{T^{(i)}}(v_i) f(v_i)-f(v_{i-2})- f(v_{i+2})-\sum_{w \in N_{T^{(i)}}(v_i)\setminus\{v_{i-2}, v_{i+2}\}}f(w).
\end{align*}
It follows that
\begin{align*}
    (d_{T^{(i)}}(v_i)-2)f(v_i) = \sum_{w \in N_{T^{(i)}}(v_i)\setminus\{v_{i-2}, v_{i+2}\}}f(w).
\end{align*}
By~\cref{cla:f_decreases}, $f(v_i)> f(v_{i+2}) \geq f(w)$ for all $w \in N_{T^{(i)}}(v_i)\setminus\{v_{i-2}, v_{i+2}\}$, which contradicts the above equation.


In conclusion, we cannot perform any Switching, Jumping, or Shifting operation for $T'$. 
Therefore, $T'$ must be isomorphic to $T^*$.
\end{case}

\begin{case}
    $t = m \geq 2$.
\begin{enumerate}
\item $t = b$
                
Since $m = b = t$, $T$ is a tree with $m$ interior vertices, $m$ leaves, and matching number $m$. 
As $n=2m$, it is clear that every interior vertex is adjacent to exactly one leaf. 
Due to~\cref{tree-up-bound}, we have $\lambda_1(T) = 1$, which implies that all the trees satisfying that every interior vertex is adjacent to exactly one leaf in this case are the extremal trees.

\item $t < b$
                
The proof is similar to~\cref{case1:m_b_extremal tree}.
\end{enumerate}
\end{case}

\begin{case} 
$t=1$
 
If $t = 1$, then $n = 2m + b - 1$. 
Through the operations, we similarly obtain a sequence of trees $$T^{(1)}\to T^{(2)} \to T^{(3)} \to \cdots \to T^{(k-2)}$$ such that  $d_{T^{(k-2)}}(v_i) = 2$ for any $1 \leq i \leq k-2$, $\card{\mathcal{V}_{T^{(k-2)}}}= 2$, and
\begin{align*}
    R_{T^{(k-2)}}(f) \leq \cdots \leq R_{T^{(2)}}(f) \leq R_{T^{(1)}}(f).
\end{align*}
Now, we perform operations on $T^{(k-2)}$.  
If $\card{N_{T^{(k-2)}}(v_{k-1})\cap B} = 1$, then there is nothing to do, and let $T^{(k-1)}$ denote $T^{(k-2)}$.
Otherwise, $\card{N_{T^{(k-2)}}(v_{k-1})\cap B} \geq 2$, and we denote these leaves by $v_{k-1,1}, v_{k-1,2}, \cdots$.
Shift these leaves $v_{k-1,2}, \cdots$ to $v_k$, and obtain a new tree $T^{(k-1)}$ with $\card{\mathcal{V}_{T^{(k-1)}}}=2$.  
Since $f(v_{k-1}) \geq f(v_k) > f(v_{k-1,2}), \cdots$, by~\cref{lem:Shifting}, we have
\begin{align*}
    R_{T^{(k-1)}}(f) \leq R_{T^{(k-2)}}(f).
\end{align*}
Note that $T(2m -3,2,b) \cong T^{(k-1)}$.
The uniqueness can be proved similarly to~\cref{case1:m_b_extremal tree}.
\end{case}
\end{proof}
\endtheorem

\begin{remark}
    The operations in our proof do not guarantee the matching number of trees to be fixed, whereas both of the original tree $T$ and the final tree $T^{*}$ can be ensured to have matching number $m$.
\end{remark}

\section{Concluding remarks}\label{section:Concluding remark}
First, we show that~\cref{Klobürštel theorem} is a corollary of~\cref{thm:m_b_extremal tree}.
\begin{proof}[\bf{{Proof of~\cref{Klobürštel theorem}}}]
According to~\cref{thm:m_b_extremal tree}, it suffices to compare eigenvalues of all the extremal graphs with different $t$, where $t=2m+b-n$.
For $t=1,2$, $T(n-2-b, 2, b)$ is the unique tree that has the Faber-Krahn property in $\mathcal{T}(n,m,b)$. 
Clearly, $P_{n-b+1} \subseteq T(n-2-b, 2, b)$.
By~\cref{thm:subgraph_larger_original graph}, we have
    \begin{align*}
        \lambda_1(T(n-2-b, 2, b))\leq \lambda_1(P_{n-b+1}) = 2(1-\cos\frac{\pi}{n-b})<1.
    \end{align*}
The first inequality is strict by~\cref{lem:simple eigenvalue-positive}.
For $t = m = b$, we have $\lambda_1(T)=1>\lambda_1(T(n-2-b, 2, b))$.
For  $3 \leq t \leq \min\{b,m\}$, 
by~\cref{thm:lower_bound_2}, we have
\begin{align*}
    \lambda_1(T(n-t-b, t, b)) \geq 2(1 - \cos\frac{\pi}{n-t-b+3}).
\end{align*}
It follows that $\lambda_1(T(n-2-b, 2, b)) < \lambda_1(T(n-t-b, t, b))$. 
Therefore, $T(n-2-b, 2, b)$ is the unique tree that has the Faber-Krahn property in $\mathcal{T}(n,b)$, which is a comet. 
\end{proof}

Lin, Liu, You and Zhao~\cite{lin2025estimatesdirichleteigenvaluegraphs}
showed that if $T = (V, E)$ is a tree with leaves as boundary and $r$ is the inscribed radius of $T$, then
    \begin{align*}
        \lambda_1(T) \geq \frac{1}{r^2}.
    \end{align*}
This implies a lower bound for a tree with diameter $D$, that is $\lambda_1(T) \geq \frac{4}{D^2}$, which only relies on the parameter $D$. 
It is natural for us to consider a more precise bound based on diameter.
This leads to the Faber-Krahn type problem for trees with given diameter $D$.
Let \begin{align*}
    &\mathcal{T}_n^D = \{T: T \text{ is a tree with } n \
    \text{vertices}\text{ and diameter } D\}.
\end{align*}
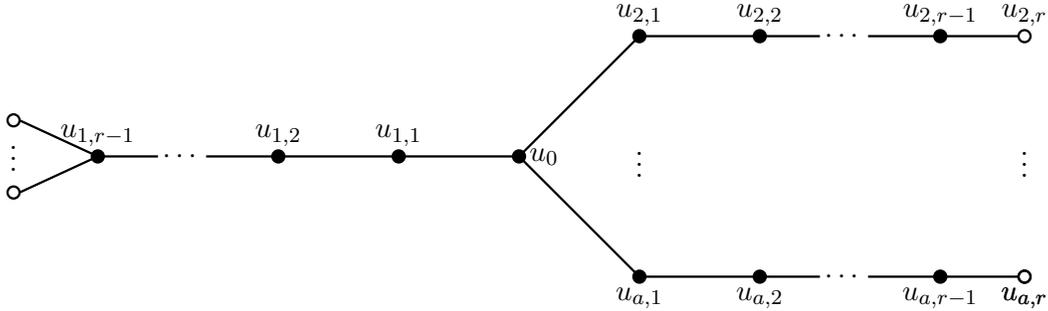
\begin{figure}[htbp]
    \centering
    \begin{tikzpicture}[x=2.0mm, y=2.0mm, inner xsep=0pt, inner ysep=0pt, outer xsep=0pt, outer ysep=0pt, scale=0.8]
    \definecolor{L}{rgb}{0,0,0}    
\definecolor{F}{rgb}{0,0,0}
\path[line width=0.30mm, draw=L, fill=F] (20,65) circle (1 mm) node[right=4pt] {$u_0$};
\path[line width=0.30mm, draw=L, fill=F] (30,75) circle (1 mm) node[above=4pt] {$u_{2,1}$};    
\path[line width=0.30mm, draw=L, fill=F] (40,75) circle (1 mm) node[above=4pt] {$u_{2,2}$};     
\path[line width=0.30mm, draw=L, fill=F] (55,75) circle (1 mm) node[above=4pt] {$u_{2,r-1}$};    
\path[line width=0.30mm, draw=L] (62,75) circle(1 mm) node[above=4pt] {$u_{2,r}$};    
\draw[line width=0.30mm, draw=L] (20,65) -- (30,75);    
\draw[line width=0.30mm, draw=L] (30,75) -- (45,75);    
\draw[line width=0.30mm, draw=L] (49,75) -- (55,75);    
\draw[line width=0.30mm, draw=L] (61.5,75) -- (55,75);    
\node at (47,75) {$\cdots$}; 
\path[line width=0.30mm, draw=L, fill=F] (10,65) circle (1 mm) node[above=4pt] {$u_{1,1}$};    
\path[line width=0.30mm, draw=L, fill=F] (0,65) circle (1 mm) node[above=4pt] {$u_{1,2}$};    
\path[line width=0.30mm, draw=L, fill=F] (-15,65) circle (1 mm) node[above=4pt] {$u_{1,r-1}$};    
\path[line width=0.30mm, draw=L] (-22,68) circle (1 mm);    
\path[line width=0.30mm, draw=L] (-22,62) circle (1 mm);
    \draw[line width=0.30mm, draw=L] (20,65) -- (0,65);    
\draw[line width=0.30mm, draw=L] (0,65) -- (-6,65);    
\node at (-8,65) {$\cdots$};    
\node at (-22,65.5) {$\vdots$};    
\draw[line width=0.30mm, draw=L] (-15,65) -- (-10,65);    
\draw[line width=0.30mm, draw=L] (-15,65) -- (-21.5,68);    
\draw[line width=0.30mm, draw=L] (-15,65) -- (-21.5,62);    
\path[line width=0.30mm, draw=L, fill=F] (30,55) circle (1 mm) node[below=4pt] {$u_{a,1}$};    
\path[line width=0.30mm, draw=L, fill=F] (40,55) circle (1 mm) node[below=4pt] {$u_{a,2}$};    
\path[line width=0.30mm, draw=L, fill=F] (55,55) circle (1 mm) node[below=4pt] {$u_{a,r-1}$};    
\path[line width=0.30mm, draw=L] (62,55) circle (1 mm) node[below=4pt] {$u_{a,r}$};    
\node at (47,55) {$\cdots$};    
\path[line width=0.30mm, draw=L] (62,55) circle (1 mm)node[below=4pt] {$u_{a,r}$};

    \draw[line width=0.30mm, draw=L] (20,65) -- (30,55);    
\draw[line width=0.30mm, draw=L] (30,55) -- (45,55);    
\draw[line width=0.30mm, draw=L] (49,55) -- (55,55);    
\draw[line width=0.30mm, draw=L] (55,55) -- (61.5,55);

\node at (30,65) {$\vdots$};    
\node at (62,65) {$\vdots$};
\end{tikzpicture}
    \caption{A generalized fork graph $GF(a,r,n)$}
    \label{fig:generalized fork graph}
\end{figure}

\begin{definition}
    Let $a \geq 2$, $r \geq 1$ , $n \geq 3$ be integers. 
    The generalized fork graph $GF(a,r,n)$ is obtained by gluing one ends of $a$ paths together, where one path is of length $r$ and the rest $(a - 1)$ paths are of length $r$, and then attach $n-ar$ leaves to the other end of the path of length $r$. 
    See~\cref{fig:generalized fork graph}.          
\end{definition}

\begin{problem}
    A tree $T$ has the Faber-Krahn property in the class $\mathcal{T}_n^D$ 
    if and only if $T$ is either a generalized fork graph 
    $GF(\lfloor \frac{n-1}{\lfloor \frac{D}{2}\rfloor}\rfloor, \lfloor \frac{D}{2}\rfloor, n)$
    or a comet.
\end{problem}
Here we consider the cases for $D=2,3,4$.
\begin{enumerate}[leftmargin=*]
    \item  $D=2$
    \setlength{\parindent}{2em}
    
    It is clear that $T \cong K_{1,n-1}$, which can also be regarded as a comet or a generalized fork graph $GF(n-1,1,n)$.
    
    \item  $D=3$ 
     
    If $T\in \mathcal{T}_n^3$, then $T$ exactly has two interior vertices and $n-2$ leaves. 
    Therefore, by~\cref{Klobürštel theorem}, a tree $T$ has the Faber-Krahn property in the class $\mathcal{T}_n^3$ if and only if $T$ is a comet.
    
    \item  $D=4$

    Let $T \in \mathcal{T}_n^4$ be a tree  with interior $\Omega$, where $n\geq 5$.
    Let $f\in\mathbb{R}^\Omega$ be a positive eigenfunction of $\lambda_1(T)$.
    Take a path $P$ of length 4 in $T$ denoted by $v_1' \sim v_1 \sim v_0 \sim v_2 \sim v_2'$.
    Clearly, $P_4 \subseteq GF(2,2,n)$.
    By~\cref{thm:subgraph_larger_original graph}, we have
    $\lambda_1(GF(2,2,n))\leq \lambda_1(P_4) = 1$.
    The inequality is strict by~\cref{lem:simple eigenvalue-positive}.
    We claim that $v_0 \notin \mathcal{V}_T$.
    Otherwise, by~\cref{thm:lower_bound_2}, we have $\lambda_1(T) \geq 1 > \lambda_1(GF(2,2,n))$,
    a contradiction.
    Next, we choose a vertex $v^* \in \Omega$ such that $f(v^*) \leq f(v)$ for any $v \in \Omega\setminus \{v_0\}$.
    If $\card{N_T(v) \cap B} \geq 2$ for $v \neq v^*, v_0$, then we denote these leaves of $v$ by $v', v'', \cdots$, and shift $v'', \cdots$ to $v^*$.
    Finally, we obtain a generalized fork graph $GF(k-1,2,n)$.
    By~\cref{lem:Shifting}, we have $R_{GF(k-1,2,n)}(f) \leq R_{T}(f)$.
    Now we claim that $\lambda_1(GF(k-1,2,n))<\lambda_1(T)$. 
    If $\lambda_1(GF(k-1,2,n))=\lambda_1(T)$, then $f$ is the common eigenfunction by~\cref{lem:extreme_eigenfunction}.
    Let $v^{**} \in \mathcal{V}_T \setminus\{v^*\}$ be a vertex whose leaves are shifted away. 
    Clearly, $N(v^{**}) \cap \Omega = \{v_0\}$.
    By the definition of eigenfunction, we have
    \begin{align*}
        \lambda_1(GF(k-1,2,n)) f(v^{**}) = d_{GF(k-1,2,n)}f(v^{**}) - f(v_0) = \lambda_1(T) f(v^{**}) = d_Tf(v^{**}) - f(v_0).
    \end{align*}
    Hence, $d_{GF(k-1,2,n)}(v^{**}) = d_T(v^{**})$, which contradicts the choice of $v^{**}$.
    Therefore, we only consider the extremal graph among $GF(k-1,2,n)$.

    The Dirichlet Laplacian matrix associated with $GF(k-1,2,n)$ is shown as follows.
    \begin{align*}
        D_k=
        \begin{pmatrix}
        k-1 & -1 & -1 & \cdots & -1\\
        -1  & n+3-2k & 0 & \cdots & 0\\
        -1  & 0 & 2 & \cdots & 0\\
        \vdots & \vdots & \vdots & \ddots & \vdots\\
        -1  & 0 & 0 & \cdots & 2
        \end{pmatrix}.
    \end{align*}   
    Its characteristic polynomial is
    $(\lambda - 2)^{k-3} P(\lambda, k)$, where
    \begin{align*}
        P(\lambda, k)
        = & \lambda^{3} + (k - n - 4)\lambda^{2} + (-2k^{2} + kn + 2k + n + 2)\lambda + (2k^{2} - kn - 3k + 2).
    \end{align*}
    Since $P(0, k)= 2k^{2} - kn - 3k + 2 \leq 2k^2 - (2k - 1)k - 3k + 2 = -2k + 2 < 0$ ($n \geq 2k - 1$ and $k \geq 3$), $P(1, k)= 1>0$.
    we have $0<\lambda_1(GF(k-1,2,n))<1$. 
    A direct computation shows that
    \begin{align*}
        P(\lambda, k+1)-P(\lambda, k) = (\lambda - 1)(\lambda - (4k - n - 1)).
    \end{align*} 
    If $P(\lambda, k+1)-P(\lambda, k) > 0$, then $\lambda_1(GF(k,2,n))< \lambda_1(GF(k-1,2,n))$. 
    If $P(\lambda, k+1)-P(\lambda, k) < 0$, then $\lambda_1(GF(k,2,n))> \lambda_1(GF(k-1,2,n))$.
    Since $\lfloor \frac{n+1}{2}\rfloor \geq k \geq 3$, it suffices to consider the cases $k = 3$ and $k = \lfloor \frac{n+1}{2}\rfloor$.
    We compute the expression of $P(\lambda, k)$ for three different values of $k$ as follows.
    \begin{align*}
        &P(\lambda, 3)= \lambda^3 - (n+1)\lambda^2 + (4n-10)\lambda - 3n + 11,\\
        &P(\lambda, \frac{n+1}{2}) = \lambda^3 - \frac{n+7}{2}\lambda^2 + \frac{3n+5}{2}\lambda - n + 1,\\
        &P(\lambda, \frac{n}{2})= \lambda^3 - \frac{n+8}{2}\lambda^2 + 2(n+1)\lambda - \frac{3n}{2} + 2.
    \end{align*}
When $n=5$ or $6$, the generalized fork graph can only be $GF(2,2,n)$.
When $n \geq 7$ is odd, $P(\lambda, \frac{n+1}{2})-P(\lambda, 3)= \frac{1}{2}(\lambda-1)(\lambda-4)(n-5)>0$.
When $n \geq 8$ is even, $P(\lambda, \frac{n}{2})-P(\lambda, 3) = \frac{1}{2}(n-6)(\lambda-1)(\lambda-3)>0$.
Therefore, we have $\lambda_1(GF(\lfloor \frac{n-1}{2}\rfloor,2,n)) < \lambda_1(GF(k-1,2,n))$ for $3 \leq k \leq \lfloor \frac{n+1}{2}\rfloor-1$.

In conclusion, a tree $T$ has the Faber-Krahn property in the class $\mathcal{T}_n^4$ 
if and only if $T$ is the generalized fork graph $GF(\lfloor \frac{n-1}{2}\rfloor, 2,n)$.
\end{enumerate}

A natural generalization of trees is planar graphs.  We know that if $G$ and $G+e$ have the same boundary vertices, then $\lambda_1(G+e,B)\geq \lambda_1(G,B)$.
So the Faber-Krahn property for planar graphs is the same as trees with the same boundary vertices.
A maximal planar graph is a planar graph with size $3n-6$. So we have the following problem.
\begin{problem}
Characterize the Faber-Krahn property for maximal planar graphs with boundary.
\end{problem}

\section*{Declaration of competing interest}

The authors declare that they have no known competing financial interests or personal relationships that could have appeared to influence the work reported in this paper.

\section*{Data availability}

No data was used for the research described in the article.

\section*{Acknowledgements}
Huiqiu LIN was supported by the National Natural Science Foundation of China (No. 12271162, No. 12326372), and Natural Science Foundation of Shanghai (No. 22ZR1416300 and No. 23JC1401500) and The Program for Professor of Special Appointment (Eastern Scholar) at Shanghai Institutions of Higher Learning (No. TP2022031).


\end{document}